\documentstyle{amsppt}
\pagewidth{5in}
\pageheight{7.8in}
\magnification=\magstep1
\hyphenation{co-deter-min-ant co-deter-min-ants pa-ra-met-rised
pre-print fel-low-ship}
\def\leaderfill{\leaders\hbox to 1em{\hss.\hss}\hfill}

\def\I{{\Cal I}}

\def\idest{i.e.,\ }
\def\a{{\alpha}}

\def\d{{\delta}}

\def\l{{\lambda}}

\def\t{{\tau}}

\def\sp{\text{\rm span}}

\def\extln{{D_n}}
\def\dnq{{D_n(q)}}
\def\dnp{{D_n^+}}
\def\dnph{{\widehat{D}_n^+}}
\def\dnj{{D_n[J]}}
\def\itq{{I_t(q)/I_{t-2}(q)}}
\def\tln{{TL(\widehat{A}_{n-1})}}
\def\otln{{O_n}}
\def\annn{\text{\rm Ann({\bf n})}}
\def\boxit#1{\vbox{\hrule\hbox{\vrule \kern3pt
\vbox{\kern3pt\hbox{#1}\kern3pt}\kern3pt\vrule}\hrule}}
\def\rabbit{\vbox{\hbox{\kern0pt
\vbox{\kern0pt{\hbox{---}}\kern3.5pt}}}}

\def\Ext{\text{\rm Ext}}
\def\Hom{\text{\rm Hom}}

\def\tableau#1{
        \hbox {
                \hskip -10pt plus0pt minus0pt
                \raise\baselineskip\hbox{
                \offinterlineskip
                \hbox{#1}}
                \hskip0.25em
        }
}

\def\tabCol#1{
\hbox{\vtop{\hrule
\halign{\strut\vrule\hskip0.5em##\hskip0.5em\hfill\vrule\cr\lower0pt
\hbox\bgroup$#1$\egroup \cr}
\hrule
} } \hskip -10.5pt plus0pt minus0pt}

\def\CR{
        $\egroup\cr
        \noalign{\hrule}
        \lower0pt\hbox\bgroup$
}



\topmatter
\title On Representations of Affine Temperley--Lieb Algebras, II
\endtitle

\author K. Erdmann and R.M. Green \endauthor
\affil
Mathematical Institute\\ Oxford University\\ 24--29 St. Giles'\\
Oxford OX1 3LB\\ England\\ 
{\it E-mail:} erdmann\@maths.ox.ac.uk
Department of Mathematics and Statistics\\ Lancaster University\\
Lancaster LA1 4YF\\ England\\
{\it  E-mail:} r.m.green\@lancaster.ac.uk
\endaffil

\abstract
We study some non-semisimple representations of affine
Temperley--Lieb algebras and related cellular algebras. 
In particular, we classify extensions between simple standard modules.
Moreover, we construct a completion which is an 
infinite dimensional cellular algebra.
\endabstract

\thanks
Part of this work was done at the Newton Institute, Cambridge,
and supported by the special semester on Representations of Algebraic Groups,
1997. The second author was supported in part by an E.P.S.R.C. postdoctoral
research assistantship.
\endthanks

\endtopmatter

\centerline{\bf To appear in the Pacific Journal of Mathematics}

\head 1. Introduction \endhead

\bigskip
The affine Temperley--Lieb algebra, $\tln$, is an infinite dimensional algebra
which occurs as a quotient of the Hecke algebra associated to a
Coxeter system of type $\widehat A_{n-1}$.  It occurs naturally in the
context of statistical mechanics \cite{{\bf 10}, {\bf 11}}, and may be thought
of as an algebra of diagrams (see \cite{{\bf 3}, \S4}).  In \cite{{\bf 6}}, the
second author used the diagram calculus and the theory of cellular
algebras as described by Graham and Lehrer \cite{{\bf 4}} to classify and
characterise most finite dimensional irreducible modules for $\tln$ and
for a larger algebra of diagrams, $\extln$. 
In fact, all simple modules for these algebras are finite dimensional;
this follows for example from \cite{{\bf 12}, 13.10.3}.

In general, a cellular algebra produces a natural class of modules, 
called {\it standard modules} or {\it cell modules}.  Special cases of 
these include Weyl modules and
Specht modules. They are not semisimple in general,
but one obtains all irreducible modules as simple quotients of
appropriate standard modules. The algebras $D_n$ (and $\tln$) have many 
finite dimensional cellular quotients.  These were used as a main tool
in \cite{{\bf 6}}, where they were called $q$-Jones algebras.

In this paper, we extend the results of \cite{{\bf 6}} by studying 
non-semisimple finite dimensional modules for $\extln$,
and also the remaining simple modules 
which were not considered in \cite{{\bf 6}}.  Moreover, we construct
completions for the algebras $\extln$ which are infinite dimensional
cellular algebras. Our techniques also work for
the algebra $\tln$, but we will not always make this explicit.
We assume $K$ is an algebraically closed field and $v \in K^*$ is such
that $\delta=v+v^{-1}$ is non-zero. 
Unless otherwise stated, we are only concerned with finite dimensional 
modules.


After reviewing some of the
results from \cite{{\bf 6}} in \S2, we classify extensions between
simple standard modules for $D_n$ in \S3.
As a consequence we can describe all finite dimensional modules
of $D_n$ which lie in blocks $B_q$ where $q$ is such that the
corresponding $q$-Jones algebra is semisimple; such blocks exist in 
abundance. In this case
every finite dimensional module in such a block is a direct sum of uniserial
modules, and an indecomposable summand has only one type of 
composition factor.

In \S4, we construct an infinite dimensional cellular algebra, $\dnp$,
which is spanned by a set of ``positive'' diagrams.  The smooth
modules for various procellular completions of this algebra (in the sense of
\cite{{\bf 7}}) provide all the uniserial modules
described in \S3, and the algebra $\extln/I_0$ embeds densely in any of
these procellular completions.  This means that the algebra $\extln/I_0$ is in
some sense almost cellular.

In \S5 we study extensions of standard modules for cellular algebras more
generally. Here we follow the approach of S. K\"{o}nig and C. Xi \cite{{\bf 9}}
where a more algebraic treatment of these algebras is given.
In \S6, we use these results to determine extensions between simple 
standard modules for the associated $q$-Jones algebras, over
arbitrary characteristic, with $\delta \neq 0$ but otherwise 
arbitrary.
Finally in \S7, we classify the remaining simple
modules for $\extln$.  

Note that our approach is to make use of the $q$-Jones algebras, which are
deformations of certain quotients of Jones' annular algebras as
introduced in \cite{{\bf 8}}.  All these algebras are finite dimensional and
cellular, with an explicit cell basis which can be dealt with by
elementary algebra \cite{{\bf 6}, \S3}.  

After submitting this paper for publication, the authors received a
copy of Graham and Lehrer's paper \cite{{\bf 5}} which derives the
results of \cite{\bf 6} independently and using different methods.
The focus of \cite{{\bf 5}} is rather different from this paper: the main
results of \cite{{\bf 5}} 
give a complete determination of the multiplicities of the composition
factors of all the cell modules, subject to some restrictions on the
ground ring.
In this paper, we are concerned with the simple standard
modules, and our central objects
of study are instead the non-semisimple modules which have filtrations by
simple standard modules.

\bigskip
\vfill\eject
\head 2. Algebras of diagrams related to affine Temperley--Lieb
algebras \endhead

In \S2 we recall from \cite{{\bf 6}}
the definitions of various algebras of diagrams which are related to
$\tln$.

\subhead 2.1 Affine $n$-diagrams \endsubhead
We assume that $K$ is an algebraically closed
field which contains a non-zero element $v$. Let $\delta = [2] = v +
v^{-1}$, we assume that $\delta $ is non-zero in $K$. We consider the
affine Temperley--Lieb algebra as a subalgebra of an algebra $\extln$ as
defined in \cite{{\bf 6}}.  This is defined in terms of diagrams as follows.

\definition{2.1.1}
An affine $n$-diagram, where $n \in {\Bbb Z}$ satisfies $n \geq 3$,
consists of two infinite horizontal rows of nodes lying at the points
$\{ {\Bbb Z} \times \{0, 1\}\}$ of ${\Bbb R} \times {\Bbb R}$,
together with certain curves, called edges, which satisfy 
the following conditions:

\item{\rm (i)}
{Every node is the endpoint of exactly one edge.}
\item{\rm (ii)}
{Any edge lies within the strip ${\Bbb R} \times [0, 1]$.}
\item{\rm (iii)}
{If an edge does not link two nodes then it is an infinite horizontal
line which does not meet any node.  Only finitely many edges are of
this type.}
\item{\rm (iv)}
{No two edges intersect each other.}
\item{\rm (v)}
{An affine $n$-diagram must be invariant under shifting to the left
or to the right by $n$.}
\enddefinition

\definition{2.1.2}
By an isotopy between diagrams, we mean one which fixes the nodes and
for which the intermediate maps are also diagrams which are shift invariant.
We will identify any two diagrams which are isotopic to each other, so
that we are only interested in the equivalence classes of affine
$n$-diagrams up to isotopy.  This has the effect that the only
information carried by edges which link two nodes is the pair of
vertices given by the endpoints of the edge.
\enddefinition

\definition{2.1.3}
Because of the condition {\rm (v)} in 2.1.1,
one can also think of affine $n$-diagrams as
diagrams on the surface of a cylinder, or within an annulus, in a natural way.
Unless otherwise specified, we
shall henceforth regard the diagrams as diagrams on the surface of
a cylinder with $n$ nodes on top and $n$ nodes on the bottom.
From now on, 
we will call the affine $n$-diagrams ``diagrams'' for short, when the
context is clear.  Under this construction, the top row of nodes
becomes a circle of $n$ nodes on one face of the cylinder, which we
will refer to as the top circle.  Similarly, the bottom circle of the
cylinder is the image of the bottom row of nodes.
\enddefinition

\example{Example 2.1.4}
An example of an affine $n$-diagram for $n = 4$ is given in Figure 1.
The dotted lines denote the periodicity, and should be identified to
regard the diagram as inscribed on a cylinder.

\topcaption{Figure 1} An affine 4-diagram\endcaption
\centerline{
\hbox to 3.638in{
\vbox to 0.888in{\vfill
        \includegraphics{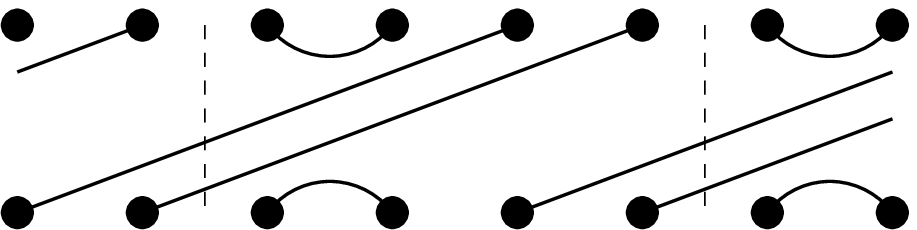}
}
\hfill}
}
\endexample

\definition{2.1.5}
An edge of the diagram $D$ is said to be vertical if it connects a
point in the top circle of the cylinder to a point in the bottom circle, and
horizontal if it connects two points in the same circle of the
cylinder.  Following \cite{{\bf 4}}, we will also call the vertical edges
``through-strings''.
\enddefinition

\definition{2.1.6}
Two diagrams, $A$ and $B$
``multiply'' in the following way, which was described in \cite{{\bf 3}, \S4.2}.
Put the cylinder for $A$ on top of the cylinder for $B$ and identify
all the points in the middle row.  This produces a certain (natural) number $x$
of loops.  Removal of these loops forms another diagram $C$ satisfying
the conditions in 2.1.1.  The product $AB$ is then defined to be $\d^x C$.
It is clear that this defines an associative multiplication.
\enddefinition

\definition{2.1.7}
Let $R = R'[v, v^{-1}]$ be the ring of Laurent polynomials over an
integral domain.
We define the associative 
algebra $\extln$ over $R$ to be the $R$-linear span of all the affine 
$n$-diagrams, with multiplication given as above.
Similarly, we may define $\extln$ over the algebraic closure of the
field of fractions of $R$, or more generally, over $K$.
\enddefinition

\subhead 2.2 Generators and relations \endsubhead

The diagram $u$ defined below plays an important r\^ole in
classifying the modules we are interested in.

\definition{2.2.1}
Denote by $\bar{i}$ the congruence class of $i$ modulo $n$, taken from
the set $\text{\bf n} := \{1, 2, \ldots, n\}$.  We index the nodes in
the top and bottom circles of each cylinder by these congruence classes in the
obvious way.
\enddefinition

\definition{2.2.2}
The diagram $u$ of $\extln$ is the one satisfying the property that for
all $j \in \text{\bf n}$, the point $j$ in the bottom circle is connected
to point $\overline{j+1}$ in the top circle by a vertical edge taking
the shortest possible route.
\enddefinition

In the case $n = 4$, the element $u$ is as shown in Figure 2.

\topcaption{Figure 2} The element $u$ for $n = 4$ \endcaption
\centerline{
\hbox to 3.638in{
\vbox to 0.888in{\vfill
        \includegraphics{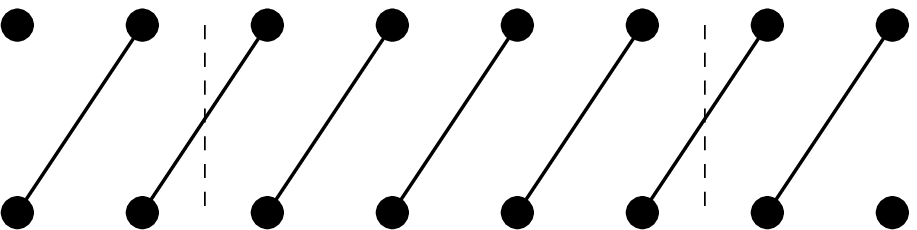}
}
\hfill}
}

\definition{2.2.3}
The diagram $E_i$ (where $1 \leq i \leq n$) has
a horizontal edge of minimal length connecting $\bar{i}$ and 
$\overline{i+1}$ in each of the circles of the cylinder, and a
vertical edge connecting $\bar{j}$ in the top circle to $\bar{j}$ in
the bottom circle whenever $\bar{j} \ne \bar{i}, \overline{i+1}$.
\enddefinition

A typical cylindrical representation of a diagram $E_i$ is shown in
Figure 3.

\topcaption{Figure 3} The element $E_2$ for $n = 5$ \endcaption
\centerline{
\hbox to 1.652in{
\vbox to 1.805in{\vfill
        \includegraphics{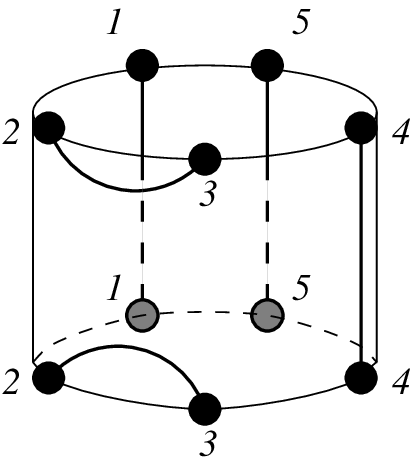}
}
\hfill}
}

\proclaim{Proposition 2.2.4}
The algebra $\extln$ is generated by elements $$E_1, \ldots, E_n, u,
u^{-1}.$$  It is subject to the following defining 
relations: $$\eqalignno{
E_i^2 &= \d E_i, & (1)\cr
E_i E_j &= E_j E_i, \quad \text{ if $\bar{i} \ne \overline{j \pm
1}$}, & (2) \cr
E_i E_{\overline{i\pm 1}} E_i &= E_i. & (3)\cr
u E_i u^{-1} &= E_{\overline{i+1}}, & (4)\cr
(u E_1)^{n-1} &= u^n . (u E_1). & (5)\cr
}$$
\endproclaim

\demo{Proof}
This is \cite{{\bf 6}, Proposition 2.3.7}.
\qed\enddemo

The affine Temperley--Lieb algebra embeds in the algebra $\extln$ as a
subalgebra.

\proclaim{Proposition 2.2.5}
The algebra $\tln$ is the subalgebra of $\extln$ spanned by diagrams, $D$,
with the following additional properties: 

\item{\rm (i)}
{If $D$ has no horizontal edges, then $D$ is the identity diagram, in
which point $j$ in the top circle of the cylinder is connected to point
$j$ in the bottom circle for all $j$.}
\item{\rm (ii)}
{If $D$ has at least one horizontal edge, then the number of
intersections of $D$ with the line $x = i + 1/2$ for any integer $i$
is an even number.}

Equivalently, $\tln$ is the unital subalgebra of $\extln$ generated by
the elements $E_i$ and subject to relations (1)--(3) of Proposition 2.2.4.
\endproclaim

\demo{Proof}
See \cite{{\bf 6}, Definition 2.2.1, Proposition 2.2.3}.
\qed\enddemo

One would like to understand the finite dimensional simple modules for
the affine Temperley--Lieb algebras $\tln$.  These are more or less
the same as the finite dimensional modules for $\extln$, as is
explained in detail in \cite{{\bf 6}, \S4}.  The above relations show
that the element $u^n$ is central in the algebra $\extln$, so by
Schur's Lemma, it acts as a scalar on any finite dimensional simple
module, and the scalar must be invertible in $K$.

The finite dimensional simple $\extln$-modules over $K$ are precisely the
simple modules for quotient algebras $\extln/\langle u^n - q \rangle$
where $q \in K$ is invertible.  These are essentially the $q$-Jones algebras of
\cite{\bf 6}, and are finite dimensional.  In \cite{\bf 6}, it was
shown that the $q$-Jones algebras are cellular, which allowed one to
classify the irreducible representations; we will summarize these
results below in \S2.5.  In this paper, we will finish this
classification and study more general finite dimensional modules.

For our approach, the finite dimensional $q$-Jones algebras are the
starting point.  The representation theory of these algebras is in
some sense similar to the representation theory of finite dimensional
quotients of the polynomial ring $K[X]$, and may be treated with
elementary methods.

\subhead 2.3 The oriented subalgebra $\otln$ \endsubhead

We also introduce a subalgebra $\otln$ of $\extln$, called the
``oriented subalgebra''.  This is related
to $\extln$ in the same way as the alternating groups are related to
the symmetric groups.  The algebra $\otln$ is defined in terms of
number of intersections with certain vertical lines; for this purpose,
we do not count a curve tangent to a line as an intersection.

\definition{2.3.1}
Define $\otln$ to be the $K$-submodule of $\extln$ spanned by diagrams $D$
such that the number of
intersections of $D$ with the line $x = i + 1/2$ for any integer $i$
is an even number.

Such diagrams are said to have the even intersection property.
Conversely, if all the numbers of intersections of a diagram $D$ with
the lines $x + 1/2$ are odd, $D$ is said to have the odd intersection property.
\enddefinition

\proclaim{Lemma 2.3.2}
Any diagram $D$ of $\extln$ has the odd intersection property or the
even intersection property.

Furthermore, if $D$ has the even (respectively, odd)
intersection property, then $u.D$ has the odd (respectively, even)
intersection property.
\endproclaim

\demo{Proof}
This was shown in \cite{{\bf 6}, Proposition 2.3.4}.
\qed\enddemo

\proclaim{Lemma 2.3.3}
The $K$-module $\otln$ is a subalgebra of $\extln$, and is generated
by $\tln$ and the elements $\{u^{2m} : m \in {\Bbb Z}\}$.
\endproclaim

\demo{Proof}
This is \cite{{\bf 6}, Lemma 2.3.3}.
\qed\enddemo

\subhead 2.4 Annular involutions \endsubhead

Let $D$ be an affine $n$-diagram associated to
the algebra $\extln$.  Throughout \S2.4
we are only concerned with diagrams $D$ with $t > 0$ vertical
edges.  (We deal with the case $t = 0$ in \S7.)
If $t > 0$, we can define the winding number $w(D)$ as follows.

\definition{2.4.1}
Let $D$ be as above.  
Let $w_1(D)$ be the number of pairs $(i, j) \in {\Bbb Z} \times {\Bbb Z}$
where $i > j$ and $\overline{j}$ in the bottom circle of
$D$ is joined to $\overline{i}$ in the top circle of $D$ by an edge
which crosses
the ``seam'' $x = 1/2$.  We then define $w_2(D)$ similarly but
with the condition that $i < j$, and we define $w(D) = w_1(D) - w_2(D)$.
\enddefinition

\remark{Remark 2.4.2}
At least one of $w_1(D)$ or $w_2(D)$ is $0$, and
$w(D)$ is always finite, although it is unbounded for a fixed value of
$n$.

The winding numbers of the diagrams in figures 1, 2 and 3 are $2$, $1$
and $0$, respectively.  Also note that $w(u^n) = n$ for any $n \in {\Bbb Z}$.

Graham and Lehrer \cite{{\bf 5}} introduce ``asymptotic modules'' $W_0$,
$W_\infty$ instead of using winding numbers.
\endremark

We recall the definition of an annular involution of the symmetric group
from \cite{{\bf 4}, Lemma 6.2}.

\definition{2.4.3}
An involution $S \in {\Cal S}_n$ is annular if and only if, for each
pair $i, j$ interchanged by $S$ $(i < j)$, we have
\item{\rm(a)}
{$S[i, j] = [i, j]$ and}
\item{\rm(b)}
{$[i, j] \cap \text{\rm Fix} S = \emptyset$ or $\text{\rm Fix}S
\subseteq [i, j]$.}

We write $S \in I(t)$ if $S$ has $t$ fixed points, and we write $S \in
\annn$ if $S$ is annular. In case $t=n$ we view the identity permutation as 
an annular involution.
\enddefinition

Using the concept of winding number, we have the following bijection.

\proclaim{Proposition 2.4.4}
Let $D$ be a diagram for $\extln$ with $t$ vertical edges ($t > 0$).
Define $S_1, S_2 \in \annn \cap I(t)$ and 
$w \in {\Bbb Z}$ as follows.

The involution $S_1$
exchanges points $i$ and $j$ if and only if $i$ is connected to $j$ in
the top circle of $D$.  Similarly $S_2$ exchanges points $i$ and $j$ if
and only if $i$ is connected to $j$ in the bottom circle of $D$.
Set $w = w(D)$.

Then this procedure produces a bijection between diagrams $D$ with at
least one vertical edge and triples $[S_1, S_2, w]$ as above.
\endproclaim

\demo{Proof}
See \cite{{\bf 6}, Lemma 3.2.4}.
\qed\enddemo

\definition{2.4.5}
For $S_1, S_2 \in \annn \cap I(t)$ where $t > 0$, 
we define $r(S_1, S_2)$ to be the smallest nonnegative integer 
satisfying $[S_1, S_2, r(S_1, S_2)] \in \otln$.
\enddefinition

\remark{Remark 2.4.6}
The integer $r(S_1, S_2)$ is necessarily $0$ or $1$ (see \cite{{\bf 6},
Definition 3.2.5}).
\endremark

\definition{2.4.7}
For $t \geq 0$ and $t \leq n$ with $ t \equiv n$(mod 2) let $I_t$ be
the span of all diagrams with at most $t$ through-strings.
This is an ideal of $D_n$. If $t > 0$ it is spanned by all elements
$[S_1, S_2, w]$ where $S_1, S_2 \in \annn \cap I(t)$. 

We also define $I_{-1} = 0$ for convenience.
\enddefinition

\definition{2.4.8 Notation} Suppose $f(X) \in K[X]$ is a polynomial, and
let $S_1, S_2 \in \annn \cap I(t)$ for $t > 0$. Then we write $$f(S_1, S_2)
$$ for the element of $I_t$ obtained from $f(X)$ by substituting 
$[S_1, S_2, k]$ for $X^k$. 
\enddefinition

\subhead 2.5 $q$-Jones algebras \endsubhead

In \cite{{\bf 6}}, finite dimensional quotient algebras $J_q(n)$ of
$\extln$ (for $n$ odd) and $\otln$ (for $n$ even) were introduced,
called $q$-Jones algebras.  These are cellular, and by using
structural results on cellular algebras, the finite dimensional
irreducible modules were classified for the case when $v$ is an
indeterminate.  Moreover, by restriction these give all
finite dimensional irreducible modules of the affine
Temperley--Lieb algebra, letting $q \in K^*$ vary.

The ordinary Jones algebra (corresponding to $q = 1$) can be simply described 
by using a homomorphism from $\extln$ to the
Brauer algebra \cite{{\bf 8}}.  There may exist a suitable deformation of
the Brauer algebra which plays a similar r\^ole for the $q$-Jones
algebra, but we do not pursue this here.

Here we allow $v$ more generally to be an element in $K^*$
(which may or may not be an indeterminate) but we assume that $\delta$
is non-zero.  In this case, the $q$-Jones algebra need not be
semisimple, but it
is still cellular. It therefore has standard modules, and all simple modules
occur as quotients of these standard modules. Therefore 
we study more generally the finite dimensional modules for $D_n$ and
their restrictions to $\tln$ which are standard modules for a $q$-Jones 
algebra, for
some $q \in K^*$. 

\definition{2.5.1} We recall the definition of the $q$-Jones algebra and 
of its standard modules. Since $u^n$ is central in $D_n$ it acts as
scalar multiplication on finite dimensional simple modules. Therefore, let
$\omega(q)$ be the ideal of $D_n$ generated by
$u^n-q$, for $q \in K^*$.  It is equal to the span of the set $$ 
\{ [S_1, S_2, w+ts] - q^s[S_1, S_2, w]: S_1, S_2 \in 
\annn\cap I(t), t \in {\Cal T}(n), w, s \in {\Bbb Z} \}.$$ Following 
\cite{{\bf 6}}, we introduce the $q$-Jones algebras as follows. 
If $n$ is odd define
$J_q(n):= D_n/\omega(q).$ On the other hand, if $n$ is even we 
define $$J_q(n) := O_n/(O_n\cap(\omega(q) + I_0))
.$$  These are finite dimensional cellular algebras, and the
restrictions of simple modules for $J_q(n)$ (viewed as $D_n$- or
$O_n$-modules) give all finite dimensional simple modules for the
affine Temperley--Lieb algebra, except simples
$M$ such that $I_0M=M$ in case $n$ is even. (We will deal with these
in \S7; see also \cite{{\bf 5}, {\bf 10}}.) 
\enddefinition

We will now describe the cellular structure and define the standard
modules.  Proofs may be found in \cite{{\bf 6}, Proposition 3.4.2,
Proposition 3.5.6}.

\definition{2.5.2} 
Assume first that $n$ is odd. Let ${\Cal T}(n) = \{ 1, 3,
\cdots, n-2, n \}$. For $t \in {\Cal T}(n)$, let
$X^t-q = \prod_{i=1}^t (X - r_i(t, q))$,
choosing some labelling for the roots.  Let $$f^{(t,j)}(X)  =\
\prod_{i>j} (X-r_i(t, q)), \ \ (1 \leq j \leq t),
$$ and set $f^{(t, 0)}(X) = 1$.
Then a cell datum for $J_q(n)$ over $K$ is given by $(\Lambda, M,
C, *)$ where
\item{(a)}{$\Lambda = \{ (t,j): t \in {\Cal T}(n), j \in ${\bf t}$
\}$, the set of weights, ordered lexicographically.}
\item{(b)}{For $(t,j) \in \Lambda$, let 
$M(t,j) = \annn \cap I(t)$.}
\item{(c)} If $S_1, S_2 \in M(t,j)$ define $$C_{S_1, S_2}^{(t,j)} := 
f^{(t,j)}(S_1, S_2^*).$$
\item{(d)}{$[S_1, S_2, w]^* = [S_2^*, S_1^*, w]$ where $S^* = w_0Sw_0^{-1}$.}
\enddefinition

\definition{2.5.3}
Assume $t > 0$ and $t\equiv n$(mod 2). The intersection 
$\omega(q) \cap O_n$ is spanned by the set
$$\{ [S_1, S_2, r(S_1, S_2) + 2w + ts] - q^s[S_1, S_2, r(S_1, S_2) + 2w]\}
$$
where $r(S_1, S_2)$ is as in 2.4.5, $S_1, S_2 \in I(t)$ and 
$w, s \in {\Bbb Z}$. 

Assume $n$ is even.  Then $J_q(n) = O_n/(O_n \cap (\omega(q) + I_0)$
is cellular. The cell data are as before but ${\Cal T}(n)$ is replaced 
by ${\Cal T}'(n) = \{ 1, 2, \cdots, n/2\}$, and in the polynomial, 
$[S_1, S_2, r(S_1, S_2) + 2k]$ is substituted for $X^k$. 

To see this, observe first that $J_q(n)$ has a basis the cosets of $$
\{[S_1, S_2, r(S_1, S_2)+ 2k]: \ \ 0 \leq k \leq t/2\},
$$ where $S_1, S_2 \in \annn\cap I(t)$. Fix
$S_1, S_2 $ and define a $K$-linear map $K[X] \to J_q(n)$ 
by $$X^k \to [S_1, S_2, r(S_1, S_2) + 2k].$$ This map takes
$X^{t/2} - q$ to $(u^n-q)[S_1, S_2, r(S_1, S_2)] \in \omega(q)$ and hence it 
induces a monomorphism $K[X]/(X^{t/2}-q) \to J_q(n)$. Then 
the set-up in 2.5.2 with $X^{t/2}-q$ instead of $X^t-q$ gives us a cell datum.
\enddefinition

\definition{2.5.4}
The filtration of $D_n$ by the ideals $I_t$ induces a filtration of
$J_q(n)$ by ideals which we also denote by $I_t$. This chain of ideals has a
refinement. Let $J_{(t,j)}$ be the span of the set $$
\{ C_{S_1, S_2}^{\lambda}: \lambda \leq (t, j), S_1, S_2 \in M(\lambda) \}.
$$ This is an ideal and we have $$
I_{t-2} \subset J_{(t,1)} \subset \ldots \subset J_{(t,t-1)}\subset
J_{(t, t)} = I_t.$$  Each quotient, as a left module, is a direct sum
of standard modules (for the general definition, see \cite{{\bf 4}} or \cite{{\bf 9}}). If $\lambda = (t, j)$ then $W(\lambda)$ can be taken as
as the span of the cosets of $\{ C_{S,T}^{\lambda}: S \in M(\lambda)
\}$ where $T \in M(\lambda)$ is fixed. We will also call it a
standard module if we view it as a module for $\extln$ or $\tln$. 
\enddefinition

\definition{2.5.5}
We note that any labelling of the roots of $X^t-q$ (respectively,
$X^{t/2} - q$) gives rise to such
chain of ideals, so in general there are many such chains. It is also
clear from \cite{{\bf 6}, \S4} that standard modules $W(t,j)$ and $W(t,k)$
are isomorphic if and only if the associated roots $r_j$ and $r_k$ are
equal.  Moreover, the standard modules depend only on the roots and
not on the particular cell chain. 
We will therefore sometimes write $W(t, \alpha)$ if $\alpha$ is a root of
$X^t-q$ (respectively, $X^{t/2} - q$). 
\enddefinition

\definition{2.5.6} 
We recall the parametrization of the simple modules
for a cellular algebra. For each $\lambda$ there is an associated
bilinear form $\phi_{\lambda}$. Let 
$\Lambda_0:= \{ \lambda: \phi_{\lambda} \neq 0 \}$. If $\lambda \in 
\Lambda_0$ then $W(\lambda)$ has a simple top, which we denote by 
$L(\lambda)$. The set $\{ L(\lambda): \lambda \in \Lambda_0 \}$ is a
full set of simple modules for a cellular algebra \cite{{\bf 4}}.

For $q$-Jones algebras in general, every standard module is
isomorphic to some $W(\lambda)$ with $\lambda \in \Lambda_0$; this is not
hard to see, and this will follow from our results in \S6. In particular
all standard modules have a simple top, and the distinct simple modules are 
in 1--1 correspondence with the distinct roots
of the polynomials $X^t-q$ for $t \in {\Cal T}(n)$. 
\enddefinition

\definition{2.5.7} For $n$ even, there is also
a cellular algebra $D_n(q) := \extln/(\omega(q) + I_0)$, with cell
structure as in 2.5.2 above where ${\Cal T}(n) = \{ 2, 4, \cdots, n
\}$.  This contains the algebra defined in 2.5.3 as a subalgebra of
half the dimension.  It is also convenient to define $D_n(q) := J_q(n)$ in
the case where $n$ is odd.
\enddefinition

\definition{2.5.8}
Similarly if $n$ is odd then there is a cellular 
finite dimensional algebra $O_n(q):= O_n/(\omega(q)\cap O_n)$. 
By considering the basis given in 2.5.2 we see that this is also cellular,
with cell datum as in 2.5.2 except that the starting polynomial is
$X^t-q^2$, and one substitutes $[S_1, S_2, r(S_1, S_2)+2k]$ into $X^k$. 
By definition, $O_n(q)$ is a subalgebra of $J_q(n)$. Both have the same
dimension, so they are isomorphic, but the isomorphism is not canonical. 
\enddefinition

\remark{Remark 2.5.9}
The algebra
$J_q(n)$ has a subalgebra 
$\Gamma_n(q)$ given as $${\Gamma \over{\omega(q) \cap \Gamma}} \text{ or }
{\Gamma \over{(\omega(q) + I_0) \cap \Gamma}}
$$ where $\Gamma = \tln$. Here the inclusion is
compatible with the cell 
structure although one must replace $I_n/I_{n-2}$ by $K$. 
\endremark

\remark{Remark 2.5.10}
Suppose $n$ is even. If $q \neq 1$ then $I_0 \subset \omega(q)$ 
automatically. But if $q=1$ then $I_0$ must be factored out. 
We will deal with this in \S7.
\endremark

\bigskip
\bigskip

\head 3. Non-semisimple modules for $\extln$. \endhead

A natural question to ask is: are there any non-trivial extensions between
simple modules for $\extln$?

There are two possible types. First there may be extensions which 
are not modules for any $\dnq$. Second, there may be extensions
which already occur for $\dnq$-modules.  We will call these
extensions of the first (respectively, second) type.
We start with the first type.
All the modules in \S3 are $\extln$-modules unless otherwise stated.
We write $\Hom(-,-)$ instead of $\Hom_A(-,-)$
and similarly for $\Ext^1$ if the algebra $A$ is clear from
the context.

\bigskip
\subhead 3.1 Some properties of extensions \endsubhead

\definition{3.1.1}
We consider finite dimensional $D_n$-modules. 
Since $u^n$ is central in $D_n$, such modules
have a natural block decomposition. For $q\in K^*$, define
$B_q$ to be the category of finite dimensional $D_n$-modules on which some
power of $u^n-q$ acts as zero, with maps all $D_n$-homomorphisms. 
That is, a finite dimensional $D_n$-module belongs
to $B_q$ if and only if $q$ is the only eigenvalue of $u^n$ on $M$. 

Thus every finite dimensional $D_n$-module $M$ is a direct sum $$
M = M_1 \oplus M_2 \oplus \cdots \oplus M_r
.$$ Here, $M_i$ belongs
to the block $B_{q_i}$ where $q_1, \ldots, q_r$ are the distinct
eigenvalues of $u^n$ on $M$. (Namely, take $M_i $ to be the kernel
of $(u^n-q_i)^{a_i}$, where $u^n$ has minimal polynomial 
$\prod_{i=1}^r (X-q_i)^{a_i}$; this is a $D_n$-module.) 

\bigskip

Suppose $M, M'$ are finite dimensional modules 
for $\extln$ which belong to blocks $B_a$ and $B_b$ respectively. 
\enddefinition

\proclaim{Lemma 3.1.2}
Suppose there is a short exact sequence $$0 \to M' \to N \to M \to 0
$$ of modules for $\extln$, where $N$ is indecomposable. Then $u^n$
has the same eigenvalue on $M$ and on $M'$.
\endproclaim

\demo{Proof}
If $a \neq b$ then $N$ would have two distinct non-zero summands
in different blocks.
\qed\enddemo

\definition{3.1.3}
Lemma 3.1.2 now allows us to restrict our attention to the case
$a=b=q$.   
\enddefinition

\proclaim{Lemma 3.1.4}
Suppose $M$ has simple top $L(\lambda)$ and assume $L(\lambda)$ does not 
occur as a composition factor of $M'$. Then there is no exact sequence
$0 \to M' \to N \to M \to 0$ such that $u^n-q$
does not act trivially on $N$.
\endproclaim

\demo{Proof}
Suppose this is false. Let $\psi$ denote multiplication by $u^n-q$.  This is
a $\extln$-endomorphism of $N$ which is non-zero, but $\psi^2$ is
zero.  Hence$$
0 \neq \text{\rm Im}(\psi) \subseteq \ker(\psi) \subseteq M'
$$ and $\psi$ induces a non-zero homomorphism from $N/M' $ to $M'$. But 
$N/M' \cong M$ has simple top $L(\lambda)$
and it follows that $L(\lambda)$ occurs in $M'$, a contradiction.
\qed\enddemo

Next, observe a more general fact.

\proclaim{Lemma 3.1.5}
Let $A$ be an algebra over some field, and let $0 \to M \to N \to M'
\to 0$ be a short exact sequence of cyclic 
$A$-modules.  

\item{(a)}{Suppose there is an ideal $I$ of $A$ such that $M = IM $
and $M' = IM'$. Then $N = IN$. 
}
\item{(b)}{Suppose $J$ is an idempotent ideal and $JM= 0$, $JM'=0$. Then
$JN=0$.}
\endproclaim

\demo{Proof} We may assume that $M \subset N$ and that $N \to M'$ is
the canonical map. Let $m$ be a generator of $M$ and $m' = n+M$ a
generator of $M'$. 

\noindent (a) Here, $M=Im$ and $M'=Im'$. Clearly, $Im + In \subseteq
N$.  Conversely, let $x \in N$.  Then
$x+M \in M' = IM' = Im'$, so there is some $z' \in I$ such that
$x+M = z'm' = z'n+M$ and then $x-z'n \in M$. Since $M = IM = Im$, there
is some $z \in I$ with $x-z'n = zm$, and hence $x \in In + Im$. 
So $N = In + Im$, and hence $N=IN$. 

\noindent (b) Let $x \in N$ and $z \in J$.  It is enough to show that
$zx=0$.  Since $J^2=J$, there
are $u_i, u_i' \in J$ such that $z = \sum_i u_iu_i'$. We have
$u_i'x +M \in JM' = 0$, so $u_i'x \in M$ for all $i$. Then $u_i(u_i'x)
\in JM = 0$ and $zx = 0.$
\qed\enddemo

\definition{3.1.6}
We concentrate now on the case when $M$ and $M'$ are standard modules
which are also simple. By Lemma 3.1.2 we may assume that $M \cong M'$, say
$M = W(t,\alpha)$. If $N$ is an extension of $M$ by itself then by
Lemma 3.1.5 we know $N = I_tN$ and $I_{t-2}N=0$. 
\enddefinition

\subhead 3.2 The map $\t_t$ \endsubhead

Lemma 3.1.5 hints that to understand the module $N$ of 3.1.6, it helps to
understand the structure of $I_t/I_{t-2}$ as a $\extln$-module.  We
therefore proceed by introducing a certain map, $\t$, on $I_t/I_{t-2}$.
The basis for $I_t/I_{t-2}$ which we work with consists of elements $D
+ I_{t-2}$ as $D$ ranges over the set of diagrams for $\extln$ with $t$
through-strings.

\definition{3.2.1}
The invertible linear map $\t = \t_t$ on $I_t/I_{t-2}$ (for $t > 0$) is
defined by its effect on the basis elements $D + I_{t-2}$ via $$
\t([S_1, S_2, w] + I_{t-2}) = [S_1, S_2, w+1] + I_{t-2}
.$$
\enddefinition

\definition{3.2.2}
The irreducible modules for $\extln$ corresponding to the section
$I_t/I_{t-2}$ may be classified using the map $\t_t$: it follows from
\cite{{\bf 6}, Proposition 4.1.1} that the irreducible
representation $\rho$ with label $(t, \a)$ satisfies $$
\rho(\t([S_1, S_2, w])) = \a\rho([S_1, S_2, w]).$$
\enddefinition

\proclaim{Lemma 3.2.3}
Consider $I_t/I_{t-2}$ as a left (respectively, right) 
$\extln$-module in the natural way.
Then $\t$ is an automorphism of left (respectively, right) \newline 
$\extln$-modules. 
\endproclaim

\demo{Proof}
It is enough to show that $\t$ is a homomorphism of
$\extln$-modules.  This is checked by considering properties of the
diagram basis.
\qed\enddemo

\definition{3.2.4}
Lemma 3.2.3 means that if $p$ is any polynomial then
$p(\t) . I_t/I_{t-2}$ is a left and right $\extln$-submodule (in
particular, an ideal) of $I_t/I_{t-2}$. 
\enddefinition

\subhead 3.3 The indecomposable modules $\I_M(S)$ \endsubhead

We now use the map $\t$ to construct explicit examples of
self-extensions of simple standard modules.

\definition{3.3.1}
Let $M$ be a simple standard  module with label $(t, \a)$.
We define the algebra $\I_M$ to be $$
{I_t/I_{t-2}}\over{(\t - \a)^2 I_t/I_{t-2}}
.$$
\enddefinition

The point of introducing $\I_M$ is that the module $N$ of 3.1.6
is naturally an $\I_M$-module.
The definition of $\t$ and Lemma 3.2.3 means that it makes sense to
apply $\t$ to $I_t/I_{t-2}$-modules, such as $M$ and $N$ (see Lemma
3.1.5).  We deduce from 3.2.2 that $(\t - \a) M = 0$.  It therefore follows
that $(\t - \a)^2 N = 0$, and thus $N$ is an $\I_M$-module.  Note that
$\I_M$ is finite dimensional.

\proclaim{Lemma 3.3.2}
The dimension of $\I_M$ is $2 (\dim M)^2$.  A basis is given by 
the images of the diagrams $$
\{[S_1, S_2, r(S_1, S_2)]\} \cup
\{[S_1, S_2, r(S_1, S_2) + 1]\}
$$ as $S_1$ and $S_2$ range over all annular involutions in the set
$I(t)$, where $M$ has label $(t, a)$ for some $a \in K^*$.
\endproclaim

\demo{Proof}
This follows from the observation that the image of 
$[S_1, S_2, w+2]$ is equal in $\I_M$ to a linear
combination of the images of $[S_1, S_2, w+1]$ and $[S_1, S_2, w]$,
because we are quotienting out by a quadratic in $\t$.
\qed\enddemo

Next we show that $\I_M$ decomposes as a direct sum of certain left
$\extln$-submodules.

\definition{3.3.3}
Let $S \in I(t)$.  We define $\I_M(S)$ to be the subspace of $\I_M$
spanned by those basis elements (as in Lemma 3.3.2) of the form $[S_1,
S, w]$ for some $S_1$ and $w$.
\enddefinition

\proclaim{Lemma 3.3.4}
As left $\extln$-modules, $$\I_M \cong \bigoplus_{S \in I(t)} \I_M (S)
.$$  Furthermore, the direct summands are all isomorphic to each
other, and each one is of dimension $2 \dim M$.
\endproclaim

\demo{Proof}
The first assertion holds because the analogous result for $I_t /
I_{t-2}$ is easily seen to be true, 
and quotienting by $(\t - \a)^2 I_t / I_{t-2}$ respects this.

To prove the second assertion we claim that for any $S, S' \in I(t)$
there exists $D \in \extln$ such that $[S_1, S, w] D = [S_1, S', w]$.
This follows from the characterisation of left cells in
\cite{{\bf 3}, Proposition 4.5.3}.  Right multiplication by the image of
$D$ in $\I_M$ establishes the required isomorphism of left
$\extln$-modules from $\I_M(S)$ to $\I_M(S')$.
\qed\enddemo

\proclaim{Lemma 3.3.5}
The modules $\I_M(S)$ are indecomposable and are isomorphic to
self-extensions of $M$.
\endproclaim

\demo{Proof}
We observe that $(\t - \a) \I_M$ is a submodule of $\I_M$ which is
isomorphic to $M$, because $(\t - \a)$ and $I_{t-2}$ annihilate it, but
$I_t$ does not, and because it has the correct dimension.  A similar
argument shows that the quotient module is also isomorphic to $M$.

If $\I_M(S)$ were decomposable, it would be isomorphic to $M \oplus
M$, and $(\t - \a)$ would annihilate it.  Such is not the case.
(Note that we use here the fact that $M$ is simple.)
\qed\enddemo

\remark{Remark 3.3.6}
There are many examples of simple standard modules.  Let $K'$ be the
algebraic closure of the prime subfield of $K$, let $v \in K$ be 
transcendental over $K'$ and
let $0 \ne \a \in K'$.  Then consideration of the cellular bilinear
form shows that the standard module with label $(t, \a)$ is simple.
\endremark

\subhead 3.4 Extensions of the first type \endsubhead

We can now classify extensions between simple standard
$\extln$-modules in the case where there are no nontrivial extensions
between these modules as $\dnq$-modules.  (The other cases will be
examined in \S6.)

\proclaim{Proposition 3.4.1}
Let $M$ be a simple standard module with 
no nontrivial self-extensions as a $\dnq$-module.
Let $$0 \to M \to N \to  M \to 0
$$ be a short exact sequence where $N$ is indecomposable.  Then
$N$ is isomorphic to $\I_M(S)$ for some $S$, and therefore for any
$S$.  In particular, if $(u^n - a)M = 0$, then $(u^n - a)N \cong M$.
\endproclaim

\demo{Proof}
By Lemma 3.1.5 we may assume $N$ is generated as a module by some $z$.  It
is clear that $(u^n - a)^2 N = 0$.
Since $N$ is indecomposable and $(u^n - a)N \ne 0$, 
Lemma 3.3.5 shows that $N$ must be equal
to $\I_M(S) . z$ for some $S$.  Comparison of dimensions shows that
the map $x \mapsto x . z$ for $x \in \I_M(S)$ is an isomorphism of
modules, so $N$ is isomorphic to $\I_M(S)$.

The observation that $$
0 \ne (u^n - a)\I_M(S) \ne \I_M(S)
$$ proves the last assertion.
\qed\enddemo

\proclaim{Theorem 3.4.2 }
Let $M$ and $M'$ be simple standard modules, where $M$ is a
$\dnq$-module and $M'$ is a $D_n(q')$-module.

\item{(a)}{If $\Ext_{\dnq}^1(M, M) = 0$ then $\Ext_{\extln}^1(M, M) =
K$.}
\item{(b)}{Suppose $M \not\cong M'$.  If $q \ne q'$ or both 
$q = q'$ and $\Ext_{\dnq}^1(M, M') = 0$, then $\Ext_{\extln}^1(M, M') = 0$.}
\endproclaim

\demo{Proof} 
Claim (a) follows from Proposition 3.4.1.

Lemma 3.1.2 proves (b) in the case $q \ne q'$.  If $q = q'$, Lemma
3.1.4 shows that any nontrivial extension of $M$ by $M'$ is also an
extension in the category of $\dnq$-modules.  This cannot happen by
hypothesis.
\qed\enddemo

\definition{3.4.3}
Computing $\Ext_{\extln}^1(M, M')$ for arbitrary simple modules $M$
and $M'$ is much more difficult, but the results above help to reduce
the problem for $\extln$ to that of the finite dimensional algebra
$\dnq$.  We will return to the corresponding 
question for $\dnq$-modules in \S6.

We note that Graham and Lehrer \cite{{\bf 5}} have some significant
results in this direction.
\enddefinition


\subhead 3.5 Uniserial modules \endsubhead

\proclaim{Proposition 3.5.1}
Suppose $A$ is any algebra and $M$ is a finite dimensional
$A$-module.  Suppose that for all simple modules $S_i$ of $A$ which
occur as composition factors of $M$ we have $$
\Ext^1(S_i, S_j) = 0 \ \ (i \neq j),
\ \ \Ext^1(S_i, S_i) \subseteq K.
$$  Then $M$ is a direct sum of uniserial modules, and  each indecomposable
summand has only one type of composition factor.
\endproclaim

\demo{Proof}
Let $\bar{A} = A/{\frak A}$ where ${\frak A}$ is the annihilator of
$M$. Then $\bar{A}$ is finite dimensional and $M$ is also an $\bar{A}$-module.
By the hypothesis, we also have $$\Ext^1_{\bar{A}}(S_i, S_j) =0$$ and $$
\Ext^1_{\bar{A}}(S_i, S_i) \subseteq K.$$ By general theory (see
\cite{{\bf 1}}) the algebra $\bar{A}$ is Morita equivalent to a direct
sum of local uniserial 
algebras, and then every finite dimensional indecomposable module is 
uniserial with only one composition factor. By the Krull-Schmidt
theorem $M$ is a direct sum of indecomposables, hence the result.
\qed\enddemo

\proclaim{Theorem 3.5.2}
Let $A = \extln$, and let $N$ be an indecomposable $A$-module with a
filtration by simple standard modules. Suppose $N$
belongs to the block $B_q$ for $q \in K^*$. 

Then if $\dnq$ is semisimple, $N$ is uniserial with all composition
factors isomorphic to $M$.
\endproclaim

\demo{Proof} 
Lemma 3.1.2 shows that all composition factors of $N$ are $\dnq$-modules.
Since $\dnq$ is semisimple, Theorem 3.4.2 shows that Proposition 3.5.1
applies, thus proving the assertion.
\qed\enddemo

\definition{3.5.3}
These uniserial modules, which may be of arbitrary length,
can be constructed explicitly by the methods of \S3.3.  The only change is the
replacement of $(\t - \a_t)^2$ by $(\t - \a_t)^k$ for arbitrary $k$.
\enddefinition

\definition{3.5.4} Similar results hold for the algebra
$\tln$. The only difference is that $\tln$ has only one simple
module associated to the top quotient $I_n/I_{n-2}$.  In the case
where $n$ is odd, the uniserial modules for $D_n$ are also uniserial
modules for $\tln$ in the way one would expect; this is true because
it holds for the standard modules (see \cite{{\bf 6}, \S4.1}).  In the case
where $n$ is even, the r\^ole of $\t$ is played by $\t^2$, but the
techniques are the same.

Note that unlike the situation for ordinary Temperley--Lieb algebras, 
one cannot construct towers of algebras in the affine case,
because $\tln$ is not contained in $TL(\widehat{A}_{n})$ in any
obvious way.
\enddefinition

\bigskip

\head 4. The algebra $\dnp$ and its procellular completion \endhead

\subhead 4.1 The algebra $\dnp$ \endsubhead

For the purposes of \S4, it is convenient to work with a subalgebra,
$\dnp$, of $\extln$ which is spanned by certain ``positive''
diagrams.  These are defined as follows.

\definition{4.1.1}
A diagram $[S, T, w]$, where $S, T \in \annn$ have $t > 0$ fixed points
each, is said to be {\it positive} if and only if $w \geq nt$.  We
denote the span of the images of the positive diagrams in $\extln/I_0$
together with the identity by $\dnp$.
\enddefinition

\proclaim{Proposition 4.1.2}
The module $\dnp$ is a subalgebra of $\extln/I_0$.
\endproclaim

\demo{Proof}
The condition for a diagram $D$ to be positive is equivalent to the
condition that each through-string contribute at least $+n$ to the
winding number $w(D)$.  Consider a product $D_1 D_2$ in $\dnp$.  This
is equal to $\d^x D$ for some diagram $D$.  Each through-string, $e$, in $D$
arises from the concatenation of a through-string in $D_1$, an
(optional) intermediate section composed of horizontal edges from the
bottom circle of $D_1$ and the top circle of $D_2$, and a
through-string in $D_2$.  Since $D_1$ (respectively, $D_2$) has at most $n/2$
horizontal edges in the bottom (respectively, top) circle, at most $n$
edges can be involved in the optional intermediate section of $e$.
Each of these contributes either $+1$, $-1$ or $0$ to the winding number of
$D$.  Adding in the contributions from $D_1$ and $D_2$ shows that $e$
contributes at least $+n$ to  $w(D)$, as required.
\qed\enddemo

\subhead 4.2 The cellular structure of $\dnp$ \endsubhead

We now show that $\dnp$ can be given a cell datum.  (This is slightly
surprising since it is infinite dimensional.)  To do this, we first
fix $q \in K^*$.

\definition{4.2.1}
Maintain the notation of 2.5.2.
The polynomial $f^{(t, j)}_c(X)$ is defined for $t \in {\Cal T}(n)$,
$1 \leq j \leq t$ and $0 \leq c \in {\Bbb N}$ via $$
f^{(t,j)}_c(X)  =\
X^{nt} (X^t - q)^c \prod_{i>j} (X-r_i(t, q))
.$$  We order the triples $(c, t, j)$ by the condition
$(c_1, t_1, j_1) < (c_2, t_2, j_2)$ iff $$
(-c_1, t_1,  j_1) < (-c_2, t_2,  j_2)$$ in the obvious lexicographic order.
\enddefinition

\proclaim{Proposition 4.2.2}
A cell datum for $\dnp$, $(\Lambda, M, C, *)$, is given as follows.
\item{(a)}{$\Lambda$ is the poset of triples $(c, t, j)$ given in
4.2.1.}
\item{(b)}{For $(c, t,j) \in \Lambda$, let 
$M(c, t,j) = \annn \cap I(t)$.}
\item{(c)}{If $\lambda = (c,t,j)$ and $S_1, S_2 \in M(\lambda)$, 
define } $$C_{S_1, S_2}^{\l} := 
f^{(t,j)}_c(S_1, S_2^*).$$
\item{(d)}{$[S_1, S_2, w]^* = [S_2^*, S_1^*, w]$ where $S^* = w_0Sw_0^{-1}$.}
\endproclaim

\demo{Proof}
It is clear that $*$ is an anti-automorphism of $\dnp$.  The image of
$C$ is a basis for $\dnp$ because, for a fixed $t$, there is exactly one
polynomial $f^{(t, j)}_c$ of degree $d$ for $d \geq nt$.  

There is one more axiom to check, which is the following.
If $\l \in \Lambda$ and $S, T \in M(\l)$ then for all $a \in \dnp$ we have $$
a . C_{S, T}^{\l} \equiv \sum_{S' \in M(\l)} r_a (S', S) C_{S', T}^{\l}
\mod \dnp(< \l),
$$ where  $r_a (S', S)$ is independent of $T$ and $\dnp(< \l)$ is the
subspace of $\dnp$ spanned by the set $$
\{ C_{S'', T''}^{\mu} : \mu < \l, S'' \in M(\mu), T'' \in M(\mu) \}
.$$

Let $\l = (c, t, j)$.  To verify the axiom, it is enough to work in
the quotient $$
{{\dnp}\over{I_n^+(c)}}
$$
where we define $I_n^+(c):= D_n^+ \cap\langle (u^n - q)^{c+1} \rangle$. 
Note that the basis arising from $C$ is compatible with this
quotient, and $C_{S, T}^{\l}$ has nonzero image.  Also note that any
basis element in the ideal $I_n^+(c)$ is
strictly less than
$\l$ by definition of the order.  The proof follows by the observation
that the $C$-basis elements with nonzero images form a cell basis for
the quotient---this is for exactly the same reason that the $q$-Jones
algebras, which correspond to the case $c = 0$, are cellular (see 2.5.2).
\qed\enddemo

\proclaim{Lemma 4.2.3}
For each $\l \in \Lambda$, there is only a finite number of $\l'$ such
that $\l \leq \l'$.
\endproclaim

\demo{Proof}
This follows from the observation that the number of $\l = (c, t, j)$
for a fixed value of $c$ is finite.
\qed\enddemo
\vfill\eject
\subhead 4.3 The procellular algebra $\dnph$ \endsubhead

\definition{4.3.1}
A cell datum satisfying the condition in Lemma 4.2.3, where the poset
$\Lambda$ is infinite but locally finite, is said to be of
profinite type, following \cite{{\bf 7}, Definition 2.1.2}.  It was
shown in \cite{{\bf 7}, \S2.2} that in this case, it makes sense to
speak of an algebra whose elements are
sums of basis elements $C_{S, T}^{\l}$ in which infinitely
many basis elements may occur with nonzero coefficients, and where the
multiplication carries over in the obvious way.  This is known as
the procellular completion, and is formally defined
(following \cite{{\bf 7}, \S2.1}) as follows.
\enddefinition

\definition{4.3.2}
Let $\Lambda$ be the poset of a cell datum of profinite type for an
algebra $A$ with base ring $R$.  We denote the set of finite ideals of
$\Lambda$, ordered by inclusion, by $\Pi$.

If $P \in \Pi$, we write $A_P$ for the cellular quotient of $A$ with basis
parametrised by the set $$
\{ C_P(S, T) : S, T \in M(\l), \l \in P \}
.$$  We write $I_P$ for the $R$-submodule of $A$ spanned by
all elements $C(S, T)$ where $S, T \in M(\l)$ for $\l
\not\in P$.  

The procellular completion, $\widehat{A}$, is the inverse limit of the
$A_P$ with the obvious homomorphisms $A_{P_1} \rightarrow A_{P_2}$
whenever $P_1 \supset P_2$.
\enddefinition

\proclaim{Proposition 4.3.3}
The procellular completion $\dnph$ of $\dnp$ is canonically \newline
isomorphic to the inverse limit $$
\cdots
\twoheadrightarrow
{{\dnp}\over {I_n^+(c)}}
\twoheadrightarrow
{{\dnp}\over{I_n^+(c-1)}}
\twoheadrightarrow
\cdots
\twoheadrightarrow
{{\dnp}\over{I_n^+(0)}}
.$$
\endproclaim

\demo{Proof}
As mentioned earlier, we may identify elements of $\dnph$ with sums $$
\sum r_{S, T}^{\l} C_{S, T}^{\l}
,$$ where in general infinitely many of the $r_{S, T}^{\l}$
are nonzero.  Because the $C$-basis is compatible with the ideal
chain given (see the proof of Proposition 4.2.2), the isomorphism
works as claimed.
\qed\enddemo

\subhead 4.4 Dense subalgebras \endsubhead

\definition{4.4.1}
It was shown in \cite{{\bf 7}, \S2} that a base of open neighbourhoods
of $0$ in a procellular algebra $\widehat A$ 
is given by the ideals $\widehat
I_P$, where $P$ is an ideal of $\Lambda$ and an element of $\widehat
A$ is in $\widehat I_P$ if and only if it is of the form $$
\sum r_{S, T}^{\l} C_{S, T}^{\l}
$$ where $r_{S, T}^{\l} \ne 0 \Rightarrow \l \not\in P$.
\enddefinition

\definition{4.4.2}
With respect to this topology, $\widehat A$ is a complete Hausdorff
topological ring with a homeomorphism
$\widehat *$ induced by the maps $*$ on the finite dimensional
cellular quotients \cite{{\bf 7}, Proposition 2.2.4}.  Furthermore,
the original cellular algebra $A$ occurs as a dense subalgebra in the
obvious way.
\enddefinition

In the case of $\dnp$, one can say more.  The following result shows
that $\extln/I_0$ is ``almost cellular''.

\proclaim{Proposition 4.4.3}
The obvious map from $\extln$ to $\dnph$ identifies the quotient
$\extln/I_0$ with a dense subalgebra of $\dnph$.
\endproclaim

\demo{Proof}
The denseness assertion is immediate from 4.4.2 and the fact that
$\dnp$ is a subalgebra of $\extln/I_0$.  We need to prove the
faithfulness. 
Consider $f \in D_n$ which involves only $[S_1, S_2, m]$, $m \in \Bbb{Z}$,
for fixed $S_1, S_2 \in \annn \cap I_t$ ($t > 0)$. It
suffices to show that any such $f$ maps to a non-zero element in 
$\dnph$. Choose $k$ such that $|m| < tk$ for all $m$ such that
$[S_1, S_2, m]$ occurs in $f$. Then these $[S_1, S_2, m]$ are linearly 
independent modulo $\langle (u^n-q)^{2k}\rangle$.
Observe that the elements $[S_1, S_2, a + b]$, for fixed $a$ and $0
\leq b < 2kt$, are independent in $$ 
{{\extln} \over {\langle (u^n - q)^{2k} \rangle + I_0}}
,$$  and hence in $$
{{\dnp} \over {I_n^+(2k-1)}}
,$$ and in the inverse limit of Proposition 4.3.3.  The claim now follows.
\qed\enddemo

\remark{Remark 4.4.4}
Notice that there are many completions of $\dnp$, since the value of
$q$ may be altered, but Proposition 4.4.3 holds in each case.
\endremark

\proclaim{Proposition 4.4.5}
Let $K_m$ be the ideal of $\dnph$ generated topologically by 
\newline $I_n^+(m-1)$. 
Then the set $\{K_m: m \in {\Bbb N}\}$ is a basis of
neighbourhoods of $0$ in $\dnph$.
\endproclaim

\demo{Proof}
Because each basis element $C_{S, T}^{\l}$ for $\dnp$ fails to be in the ideal
$K_m$ for sufficiently large $m$, and the
ideals $K_m$ are spanned (topologically) by the basis elements
they contain, it follows that the set $K_m$ is a basis of
neighbourhoods for $0$ in $\dnph$.
\qed\enddemo

\definition{4.4.6}
Recall from \cite{{\bf 7}, \S2.3} that a 
smooth module for a procellular algebra is one whose annihilator is
open.  It follows from Proposition 4.4.5 that the smooth modules for
$\dnph$ are precisely those annihilated by $I_n^+(c)$ for
a sufficiently large $c$.  Thus all the uniserial modules of \S3.5
occur as smooth modules for $\dnph$ for a suitable value of $q$.
\enddefinition

\remark{Remark 4.4.7}
Finally in this section, we note that the constructions of \S4 may
be adapted to work for $\otln$ and its related algebras by making
small changes to the arguments.
\endremark

\bigskip

\head 5. Standard modules for cellular algebras
\endhead

In \S5, we take $A$ to be an arbitrary 
finite dimensional cellular algebra. The aim is to study 
extensions of standard modules. Since $A$ is general, it is convenient
to work without bases as much as possible, so we follow the approach
in \cite{{\bf 9}}.

\subhead 5.1 Cell chains \endsubhead

\definition{5.1.1}
The paper \cite{{\bf 9}} gives an equivalent definition for a cellular
algebra.  It first 
introduces  the notion of a cell ideal \cite{{\bf 9}, 3.2}.  The algebra
$A$ is then cellular if there is a chain of ideals $$\eqalignno{
0 &= J_0 \subset J_1 \subset \ldots \subset J_r = A & (*)\cr
}$$ such that $J_i/J_{i-1}$ is a cell ideal in $A/J_{i-1}$. It is proved that
a cell ideal $J$ satisfies one of the following. Either $J^2 = 0$, or else
$J = AeA$ (and $J^2=J$) where $e$ is a primitive idempotent of $A$
fixed by the involution $*$, and moreover $eAe = K$, and
multiplication in $A$ induces a bimodule isomorphism 
$J \cong Ae \otimes_K eA$. 
In the last case, $Ae$ is isomorphic to the associated standard module.
Moreover, the case $J^2=J$ is equivalent to the condition that the associated
bilinear form as in \cite{{\bf 4}} is non-zero. If this is the case, we say
that the associated standard module is ``nice''. 
\enddefinition

\remark{Remark 5.1.2}
Note also that in the second case, $J$ is a heredity ideal in the
sense of \cite{{\bf 2}}. So if all quotients are idempotent then the algebra
is quasi-hereditary.
\endremark

\definition{5.1.3} In \cite{{\bf 9}} the following is proved. Given any
cellular algebra $B$ and any $B$-module $J$ then there is at least one
cellular algebra $A$ 
which is an extension of $B$ by $J$, in which $J$ is a cell ideal with square
zero. Hence not much can be said in general for standard modules which
are not nice. 
However, the $q$-Jones algebras have a special property, namely that every
standard module is isomorphic to a nice standard module. 
\enddefinition

\subhead 5.2 Homomorphisms between standard modules \endsubhead

\definition{5.2.1}
Suppose $A$ is cellular, with cell chain (*). In this
section, we denote the standard module associated to $J_i/J_{i-1}$ by
$W(i)$.  The set of weights
is $\Lambda = \{ 1, 2, \cdots, r \}$, with the natural order,
following \cite{{\bf 4}}:
note that this is opposite to the usual order for a quasi-hereditary algebra. 
Suppose $J_i/J_{i-1}$ is idempotent, then there is a primitive idempotent
$e_i$ of $A$ such that $J_i = Ae_iA + J_{i-1}$ and moreover
$W(i) \cong Ae_i/J_{i-1}e_i$. The following is a
minor generalization of \cite{{\bf 4}, 2.6}, which is based on the
observation that
if $J$ is an idempotent ideal and $M$ is a module with $JM=0$, then for
any $a \in A$ we have $\Hom_A(Ja, M) = 0$. 
\enddefinition

\proclaim{Lemma 5.2.2}
Suppose $J_i = Ae_iA + J_{i-1}$ where $e_i$ is a primitive idempotent. 
If $i < k$, then $\Hom(Ae_i, W(k)) = 0$. In particular if $W(i)$ is
the associated standard module then $\Hom(W(i), W(k)) = 0 $.
\endproclaim

\demo{Proof}
The module $W(i)$ can be taken as a quotient of $Ae_i$, and
$W(k)$ is contained in $A/J_{k-1}$.  Hence we have $$
\Hom(W(i), W(k)) \subseteq \Hom(Ae_i, A/J_{k-1}) = 0
,$$
since $e_i$ is an idempotent, and by the observation, 
$e_i(W(k) ) \subseteq e_i(A/J_{k-1}) = 0$.
\qed\enddemo

\definition{5.2.3} This gives some information on repeated occurrences
of standard modules.
Suppose $W(i)$ is a standard module such that $J_i = Ae_iA + J_{i-1}$
for a primitive idempotent $e_i$. If $j \neq i$ and $W(j) $ is
isomorphic to $W(i)$
then by the above Lemma we have $j < i$, and moreover the associated
cell quotient $J_j/J_{j-1}$ is nilpotent. 
\enddefinition

\subhead 5.3 Extensions of standard modules \endsubhead

Suppose $W(i)$ is a standard module such that
$J_i = Ae_iA + J_{i-1}$. Consider 
$\Ext^1(W(i), W(k)).$ In case $A$ is quasi-hereditary then 
$\Ext^1(W(i), W(k))= 0$ for $i \leq k$ (see \cite{{\bf 2}}). It is natural
to ask whether this generalizes to cellular algebras. 
 
The module $W(i)$ has a projective cover
$$0 \to J_{i-1}e_i \to Ae_i \to W(i) \to 0
$$
Applying $(-, W(k)) := \Hom_A(-, W(k))$, this gives
$$0 \to (W(i), W(k)) \to (Ae_i, W(k)) \to (J_{i-1}e_i, W(k))
\to \Ext^1(W(i), W(k)) \to 0
$$
Suppose $i \leq k$. Then the first two terms are isomorphic, by
the above Lemma if $i < k$, and by \cite{{\bf 4}} for $i=k$.
Hence we get for $i \leq k$ that
$$\Ext^1(W(i), W(k)) \cong \Hom(J_{i-1}e_i, W(k))
$$
In general, this can be non-zero. 
For example, any commutative local algebra over an algebraically
closed field is cellular, with $*$ the identity, and all standard
modules are isomorphic to $K$ (see \cite{{\bf 9}, 3.5}). But
$\Ext^1(K,K) \neq 0$ unless the algebra is 1-dimensional.

\bigskip 
\subhead 5.4 Idempotent cell chains \endsubhead

In particular situations one can say more. The following will be used
for the algebras $\dnq$.

\proclaim{Lemma 5.4.1}
Suppose the ideals $J_i, J_{i-a}$ in the cell chain are idempotent.
Suppose the quotient $J_i/J_{i-a}$ is the direct sum of copies of
$M_i$ where $M_i$ is generated by a primitive idempotent of $A/J_{i-a}$
and has a 
filtration with quotients isomorphic to $W(i)$. Then:

\item{(a)}{if $i < k$ then 
$\Ext^1(W(i), W(k)) = 0$;}
\item{(b)}{$\Ext^1(W(i), W(i)) \neq 0$ if and only if $M_i$ has at 
least two quotients isomorphic to $W(i)$.}
\endproclaim

\demo{Proof}
Let $U = J_{i-1}e_i.$ By \S5.3 we must look at 
$\Hom(U, W(k))$ for $i \leq k$. 
By the hypothesis, $M_i$ is the projective cover of $W(i)$ as a module
over $A/J_{i-a}$, and there is an exact sequence $$
0 \to V \to M_i \to W(i) \to 0
,$$ where $V$ has a filtration with quotients isomorphic to $W(i)$. The 
projective cover over $A$ (from 5.3)  factors through this, and we get $$
0 \to J_{i-a}e_i \to U \to V \to 0
$$ 
and an exact sequence $$
0 \to \Hom(V, W(k)) \to \Hom(U, W(k)) \to \Hom(J_{i-a}e_i, W(k)).
$$ Since $J_{i-a}$ is idempotent and $J_{i-a}W(k) = 0$ for $i \leq k$
we have
$\Hom(J_{i-a}e_i, W(k)) = 0.$

\noindent (a) Assume $i < k$.  Then $\Hom(W(i), W(k)) = 0$ (see Lemma 5.2.2),
and since $V$ is filtered with quotients $W(i)$, it follows by induction
that $\Hom(V, W(k)) = 0$ and hence that $\Hom(U, W(k)) = \Ext^1(W(i), W(k))=0$.

\noindent (b) Let $i=k$. If $V=0$
then the extensions are
zero. Otherwise, 
there is a surjective map from $V$ onto $W(i)$ and hence $\Ext^1(W(i),
W(i))$ is non-zero. 
\qed\enddemo

\bigskip

\head 6. Extensions in the non-separable case \endhead

We will now study standard modules for the algebras $A :=\dnq$
by modifying the cell chain described in 2.5.4 and 2.5.7. The same results
hold with minor modifications also for the algebras
$O_n(q)$, $\Gamma_n(q)$ and also for 
$D_n^+/I_n^+(c)$. To keep the notation
simpler we work with $D_n(q)$; for the modifications, see 6.3.4, 
6.3.5 and 6.3.7.  We
are interested in the structure of the section $\itq$,
where $I_t(q)$ is the ideal of $A$ corresponding to the ideal $I_t$
of $\extln$ which has a good description by 2.5.1.

\subhead 6.1 Direct sum decomposition of $\itq$ \endsubhead

The structure of $\itq$ is governed by the polynomial $X^t-q$, and the
filtration given by the cell basis is an analogue of the cyclic decomposition
of $K[X]/(X^t-q)$. In order to obtain a direct sum decomposition of
$\itq$ into indecomposable direct summands, we will use the primary 
decomposition of $K[X]/(X^t-q)$, whose explicit data are well-known. 

\definition{6.1.1}
Suppose $r_1, r_2, \ldots, r_m$ are the distinct roots of $X^t-q$.  
Then we have
$X^t-q = (X^m-w)^s$ where $X^m-w = \prod_{i=1}^m (X-r_i)$.  The
ring $K[X]/(X^t-q)$ is the direct sum of $m$ local rings, namely $$
{{K[X]}\over{(X^t-q)}} \cong \bigoplus_{i=1}^m 
{{K[X]}\over{(X-r_i)^s}}.
$$ The block idempotents are the cosets of $$q_i(X) := \prod_{j \neq i}
(X-r_j)^s\cdot d_i,$$ where $d_i =  d_i(X)$ is a polynomial such that
$\sum q_i(X)=1$.
(This exists since the greatest common divisor 
of the $\prod_{j\neq i} (X-r_j)^s$
$(1 \leq i \leq m)$ is $1$, so the ideal they generate in $K[X]$ is the
whole ring.)
Moreover, for $0 \leq c \leq s-1$, let $$g_i^{(c)}(X) :=
(X-r_i)^cq_i(X).
$$ Then the $i$-th summand has $K$-basis the cosets of $g_i^{(c)}(X)$ for
$0 \leq c \leq s-1$.
Note that
we have $$
Xg_i^{(0)}(X) \equiv r_ig_i^{(0)}(X) \mod g_i^{(1)}(X).$$
\enddefinition

\definition{6.1.2} 
We will now give an analogous direct sum decomposition of $\itq$. 
For a fixed $T \in \annn\cap I(t)$, let
$U(T)$ be the space spanned by all $[S,T,k]$ where $S$ runs through the
annular involutions with $t$ fixed points, and $0 \leq k \leq
t-1$. This is a left module (modulo $I_{t-2}$) and $\itq$ is the
direct sum of the $U(T)$ where $T$ runs through the annular
involutions with t fixed points. Moreover,
all $U(T)$ are isomorphic. We now fix $T$ and study $U(T)$. 
\enddefinition

\definition{6.1.3}
Suppose $S \in \annn\cap I_t$ is fixed as well. Consider the space spanned by
all $[S, T, k]$ in $\extln$. 
We identify this with polynomials (as in 2.4.8), via
the vector space isomorphism which takes $X^k $ to $[S,T,k]$. This
maps $$
X^t-q \mapsto [S,T,t]-q[S,T,0] = (u^n-q)[S,T,0],$$ which is zero in 
$A$.
Hence the above map induces a vector space isomorphism of
$K[X]/(X^t-q)$ and the span of the $[S, T, k]$ in $A$. If $f(X)$ is a 
polynomial, write $f(S)$ for the corresponding element $f(S,T)$.
\enddefinition

\proclaim{Lemma 6.1.4} 
The space $U(T)$ has $K$-basis the cosets of $$
\{ g_i^{(c)}(S): S \in \annn\cap I(t), 1 \leq i \leq m, 0 \leq c \leq s-1 \}.
$$
\endproclaim

\demo{Proof}
First note that the number of elements equals the dimension. To
show linear independence, it suffices to consider the $g_i^{(c)}(S)$
for a fixed $S$. Under the above identification, these correspond to
the cosets of $g_i^{(c)}(X)$ in 
$K[X]/(X^t-q)$ which we know form a basis, hence the claim follows. 
\qed\enddemo

\definition{6.1.5}
We have as vector spaces $$U(T) = \bigoplus_{i=1}^m Q_i(T),$$ where
$Q_i(T)$ is the span of all $g_i^{(c)}(S)$ for fixed $i$.
\enddefinition

\subhead 6.2 The modules $Q_i(T)$ \endsubhead

It turns out that the vector space decomposition of 6.1.5 is the
decomposition of $U(T)$ into indecomposable summands.

\proclaim{Proposition 6.2.1} 
\item{(a)}{The $K$-space $Q_i(T)$ is an $A$-module and it is indecomposable.}
\item{(b)}{There is a filtration of $Q_i(T)$ given by $$
0 \subset V_{s-1} \subset \ldots \subset V_1 \subset V_0 = Q_i(T)
$$ where all the successive quotients are isomorphic.  Moreover,
for each $k$ we have $Q_i(T)/V_{s-k} \cong V_k$.
}
\endproclaim

\demo{Proof}
We first prove (a). It suffices to show that $Q_i(T)$ is invariant
under $u$ and $E_1$. First, we have $u[S, T, k] = [S^u, T, k']$
where $k' = k$ if $S$ fixes $n$, and $k' = k+1$ otherwise, and where
$S^u$ is the annular involution which takes $i$ to $S(i-1) + 1$. 
In particular this shows that $k' = k+f$ where $f = f(u,S)$ is independent 
of $k$. We deduce that $u\cdot g_i^{(c)}(S)$ is the element 
obtained from $X^f \cdot g_i^{(c)}(X)$ by substituting $[S^u, T, k]$
for $X^k$; this too lies in $Q_i(T)$. 

Consider now $E_1[S,T,k], $ that is, $[S_1, S_1, 0]\cdot
[S,T,k]$ (with $S_1 = (1 2)$). Either this lies in $I_{t-2}$ for all $k$,
 or else it is
equal to $$\delta^{(S_1, S)}[S', T, k+ f]
,$$ where $f = f(S_1, S)$ depends only on $S_1, S$ but not on $k$, and where
$S'$ is an element in $\annn \cap I_t$. Here $(S_1, S)$ is the number
of loops which appear in the product (see 2.1.6).
(To see this, we may use the map $\tau $ from \S3.2 which we know
induces a module homomorphism of $\itq$. Namely,
$[S, T, k] = \tau^k([S,T,0])$ and $\tau $ commutes with $E_1$ modulo
$I_{t-2}$.  So the number $f$ is determined by $E_1[S, T, 0]$ if this
is non-zero.)
It follows that $E_1g_i^{(c)}(S)$ is either zero or is equal to
the element obtained from $\delta^{(S_1, S)}\cdot X^fg_i^{(c)}(X)$ by 
substituting $[S', T,k]$ for $X^k$, and this again lies in $Q_i(T)$.

Now we will show that $Q_i(T)$ is indecomposable. Let 
$e := (\delta^{(T,T)})^{-1}[T,T,0]$.  Then $e$ is an idempotent which belongs
to $I_t$, and we have $$
Ae \equiv I_te \equiv U(T) \mod I_{t-2}
.$$ It follows that the endomorphism ring of $U(T)$ is isomorphic to 
$eAe$ (mod $eI_{t-2}e$). This ring is  spanned by the cosets of the $[T,T,k]$
and it is isomorphic to $${{K[X]}\over{(X^t-q)}}.$$ 
(An isomorphism is induced by $X \to (\delta^{(T,T)})^{-1}[T,T,1]$.)
We know that this has precisely $m$ primitive idempotents and hence 
$U(T)$ must be the direct sum of $m$ indecomposable modules. We have already
found $m$ non-zero summands, the $Q_i(T)$, so these must be indecomposable. 

We now tackle (b). Let $V_d$ be the span of $\{ g_i^{(c)}(S): c \geq d, 
S \in \annn \cap I(t) \}$. These are submodules which give a filtration
of $Q_i(T)$. We claim that
$V_j/V_{j+1} $ is isomorphic to $V_0/V_1$. Namely, there is a
well-defined linear map which takes 
$g_i^{(j)}(S) + V_{j+1}$ to $g_i^{(0)}(S) + V_1$.
The action of $u, E_1$ was determined above, and we see that this map
is an $A$-homomorphism, hence an isomorphism. 
Similarly we have an $A$-homomorphism $Q_i(T) \to V_k$ which takes 
$g_i^{(c)}(S)$ to $g_i^{(c+k)}(S)$ and this gives an isomorphism between
$Q_i(T)/V_{s-k}$ and $V_k$, as required. 
\qed\enddemo

\proclaim{Lemma 6.2.2}
The quotient $Q_i(T)/V_1$ is the standard module $W(\lambda)$ with label
$\lambda = (t,r_i)$.
\endproclaim

\demo{Proof}
Let $S_1, S_2 \in \annn\cap I(t)$ and 
$\omega \in Q_i(T) \setminus V_1$.  Then $$
[S_1, S_2, k+1]\omega \equiv r_i[S_1, S_2, k]\omega \mod V_1
$$
and there are $\omega, S_1, S_2$ such that this is non-zero.
This characterizes $L(\lambda)$: see 3.2.2.  (In fact, the same holds
for the standard module $W(\lambda)$.)

We may assume $\omega$ is the coset of $g_i^{(0)}(S)$ for some
$S \in \annn\cap I(t)$.  The proof then follows by the argument of (a) in
the proof of Proposition 6.2.1. If $S_1=S_2=S$ then the element is non-zero. 
\qed\enddemo

We note that 6.2.1 and 6.2.2 prove the claim made in 2.5.6 that any
standard module is isomorphic to $W(\lambda)$ for some $\lambda \in \Lambda_0$.

\remark{Remark 6.2.3}
We remark that $Q_i(T)$ is the projective cover of $W(\lambda)$ as a module for
$A/I_{t-2}$. Actually, $q_i(T)$ is a non-zero multiple of an idempotent of
$A/I_{t-2}$ which is then primitive. 
\endremark

\subhead 6.3 Properties of standard modules for $\dnq$ \endsubhead

\proclaim{Theorem 6.3.1}Let $\lambda = (t, r_j)$ and 
$\mu = (t', s_i)$  be weights. 

\item{(a)}{If $\lambda \neq \mu$ and $t \leq t'$ then 
$\Ext^1(W(\lambda), W(\mu)) = 0$.}

\item{(b)}{If $X^t-q$ is separable 
then $$\Ext^1(W(\lambda), W(\lambda)) = 0.$$  Otherwise, $$
\Ext^1(W(\lambda), W(\lambda)) = K.$$ 
}
\endproclaim

\demo{Proof}
Most of this follows from Lemma 5.4.1. The only part missing is the
dimension of $\Ext^1(W(\lambda), W(\lambda))$ if $X^t-q$ has multiple roots. 
Here we have $$\Ext^1(W(\lambda), W(\lambda)) \cong \Hom(V_1,
W(\lambda)),$$ where $V_1$ is as in Proposition 6.2.1. The
epimorphism $Q_j(T) \to V_1$ gives an inclusion 
$$
\Hom(V_1, W(\lambda)) 
\subseteq \Hom(Q_j(T), W(\lambda))
.$$  This is 1-dimensional by Remark 6.2.3
and since $L(\lambda)$ occurs only once in $W(\lambda))$. 
Hence 
we have $$
\dim \left( \Ext^1 \left( W(\lambda),
W(\lambda) \right) \right) \leq 1.
$$  It is non-zero since $Q_j(T)/V_2$ is a non-split extension (the quotient
of an indecomposable projective).
\qed\enddemo

\proclaim{Lemma 6.3.2}
If $A = \dnq$ then all standard modules are nice. Moreover, $A$ is
quasi-hereditary if and only if $A$ contains a primitive $t$-th root
of $1$ for all $t \in {\Cal T}(n)$.
\endproclaim

\demo{Proof}
By Proposition 6.2.1, all standard modules occurring in $\itq$ are
isomorphic to $W(t,i)$ for $1 \leq i \leq m$ if $X^t-q$ has $m$
distinct roots. 

The algebra $A$ is quasi-hereditary if and only if all cell quotients
are idempotent, and this the case if and only if $X^t-q$ is separable
for all $t \in {\Cal T}(n)$. 
\qed\enddemo

\remark{Remark 6.3.3}
We note also that $A$ cannot be semisimple unless all polynomials
$X^t-1$ for $t \in {\Cal T}(n)$ are separable. 
\endremark

\definition{6.3.4} 
Consider the algebra $O_n/(\omega(q) \cap O_n)$
for $O_n/((\omega(q) + I_0) \cap O_n)$. Then we have exactly the same
results for the standard modules, replacing the polynomial $X^t-q$ by
either $X^t-q^2$ if $n$ is odd, or by $X^{t/2}-q$ in case $n$ is
even. (The substitution must be made as described in 2.5.1.) The
statement in Theorem 6.3.1 remains true, provided in (b) the polynomial
is replaced accordingly. 
\enddefinition

\definition{6.3.5}
Consider the algebra $\Gamma_n(q)$, \idest the affine
Temperley--Lieb algebra modulo the ideal $\omega(q) \cap \tln$ (or 
$\omega(q) + I_0 \cap \tln$).
As we remarked in 2.5.9, this is also cellular, and it is clear that
the ideals corresponding to $I_t$ (for $t < n$) are identical with the 
ideals of the algebra in 6.3.4. 
Moreover, $I_n/I_{n-2}$ is 1-dimensional, spanned
by the coset of $1$. So for $t < n$ the standard modules are the same
as for the algebras in 6.3.4, and the only standard module for $t=n$
is the trivial module. 
\enddefinition

Let $A= J_q(n)$ and $B = \Gamma_n(q)$. The
following is easy to see, by considering the action (see Lemma 6.2.2) 
and by comparing dimensions. 

\proclaim{Lemma 6.3.6}
Suppose $\lambda = (t,a)$ is a weight for $A$. Then 
the restriction of $W(\lambda)$ to $B$ is the standard module with label
$(t, a^2)$ if $n$ is odd and $t < n$, or $(t, a)$ if $n$ is even and $t < n$.
If $t=n$ then the restriction is $K$.
\endproclaim

\bigskip

\definition{6.3.7} Similar results hold for standard modules
of the algebras $D_n^+/ I_n^+(c)$. Instead
of $I_t, I_{t-2}$ we must take the ideals 
$I_{(a, t)}$ of $D_n^+$ where
$I_{(a,t)}$ is the span of
$$ \{ 
C_{S_1, S_2}^{\lambda}: \lambda = (d, s, j), s \leq t, d \geq a \}
$$
for $a \in \Bbb{N}$ and $0 \neq t \in {\Cal T}(n)$. The non-zero quotients
induced on the factor algebra are of the form
$I_{(a,t)}/I_{(a, t-2)}$  or
$I_{(a, t_0)}/I_{(a+1, t_0 )}$ where $t_0 $ is the smallest member of
${\Cal T}(n)$. One replaces the 
polynomial $X^t-q$ by
$X^{nt}(X^t-q)^{a+1}$ or
$X^{nt_0}(X^{t_0}-a)^{a+1}$. 
\enddefinition

\bigskip

\subhead 6.4 Extensions of the second type \endsubhead

We have a surjective homomorphism from $\Gamma:= \tln$ onto $B$.  This
induces an inclusion $$
\Ext^1_B(W(\lambda), W(\mu)) \to \Ext^1_{\Gamma}(W(\lambda), W(\mu))
$$ if $W(\lambda), W(\mu)$ are standard modules for $B$, and this
subspace consists of the extensions which we called ``the second type''
in \S3.  The analogue of Theorem 6.3.1 gives
the following. 

\proclaim{Lemma 6.4.1}
Let $\lambda = (t, a)$ and $\mu = (t', b)$. Then $$
\Ext^1_B(W(\lambda), W(\mu)) = 0
$$ if $t \leq t'$ and $\lambda \neq \mu$.
Furthermore, $\Ext^1_B(W(\lambda), W(\lambda)) = K$ if $t < n$ and 
the associated polynomial has multiple roots; otherwise,
$\Ext^1_B(W(\lambda), W(\lambda)) = 0$. 
\endproclaim

This is what we can say in general, independent of the value of $\delta$. 

\bigskip

\head 7. The ideal $I_0$  for $n$ even.
\endhead

Suppose $n$ is even, so that $u^n$ acts as identity on 
the ideal $I_0$.  Since we want to study all finite dimensional simple
modules for $D_n$ on which $u^n$ acts as identity, we must
also deal with modules $M$ such that $I_0M \neq 0$.  This was not
tackled in \cite{{\bf 6}}, even in the case of simple modules.

\subhead 7.1 The structure of the ideal $I_0$ \endsubhead

We start by describing a labelling for the diagrams in $I_0$. 

\definition{7.1.1}
Let $I_n(0)$ be the set of all permutations $S$ of $\Bbb Z$ which have the 
following properties:
\item{(a)}{for all $k \in {\Bbb Z}$, $S(n+k) = S(k) + n$;}
\item{(b)}{the image of $S$ in ${\Cal S}_n$ is an annular involution
with no fixed points.}
\enddefinition

\definition{7.1.2}
With this definition, there is a one to one correspondence between the
diagrams in $I_0$ and triples $[S_1, S_2, k]$ where $S_1, S_2 \in
I_n(0)$ and where $k \geq 0$. Here $S_1, S_2$ describe the top and the
bottom of the diagram, and $k$ is the number of infinite bands.
\enddefinition

\remark{Remark 7.1.3}
The elements $I_n(0)$ can also be labelled as 
pairs $(S, x)$ where $S\in \annn\cap I(0)$ and $x$ is 
the smallest positive integer such that the seam $x+1/2$ does not
intersect any edge of the half-diagram.  Although this seems more natural,
it introduces extra notation which is not relevant, so we will
use the labelling in 7.1.2.
\endremark

\definition{7.1.4} 
We start by defining ideals of $D_n$ contained in $I_0$.
Let $g(X)$ be a polynomial in $K[X]$. For $S, T \in I_n(0)$, let
$g(S,T) \in I_0$ be the element obtained from $g(X)$ by substituting
$[S,T,k]$ for $X^k$. If $J$ is any ideal of $K[X]$ then we define a subspace
$U_J$ of $I_0$ by $$
U_J := \sp
\{ g(S,T): g(X) \in J \text{ and } S, T \in \annn\cap I(0)\}.
$$ This is an ideal of $D_n$. Consider the left action. We have
$u.g(S,T) = g(S^u,T)$, where $S^u(i) = S(i-1) + 1$ and we have $$
E_1[S,T,k] = \delta^{(S_1, S)}\cdot [S',T,k+f],
$$ for some $S' \in I_n(0)$ where $f$ does not depend on $k$. 
So $E_1 . g(S,T)$ is the element 
obtained from $$
\delta^{(S_1, S)}X^f\cdot g(X)
$$ by substituting $[S',T,k]$ for $X^k$. 
\enddefinition

\definition{7.1.5}
If $T \in \annn\cap I(0)$ is fixed, we have the idempotent
$(\delta^{(T,T)})^{-1}[T,T,0]$, and we write $$
U_J(T) := U_J[T, T, 0] = \sp \{ g(S,T): g \in J, S \in I_n(0) \}.
$$
\enddefinition

\proclaim{Lemma 7.1.6}
Suppose $M$ is a finite dimensional
simple $D_n$-module such that $I_0M \neq 0$. Then there is some ideal 
$J$ of $K[X]$ such that $U_JM = 0$. Moreover $u^n$ acts as identity 
on $M$.
\endproclaim

\demo{Proof} 
Since $M$ is simple, we have $I_0M = M$.  Note that 
$I_0$ is generated by all $[T,T,0]$ for $T \in I_n(0)$ as a left
$A$-module.  (First, for any $S \in I_n(0)$ we have
$[S,T,0]\cdot [T,T,w] = \delta^{(S,T)}\cdot [S,T,w]$. Moreover, for any $T$ 
there is some $S$ such that the group $\langle S, T \rangle$ has
exactly one infinite orbit on ${\Bbb Z}$, so that $[T,S,0]\cdot
[T,T,w]$ is a non-zero multiple of $[T,T,w+1]$.)
So there is some $[T,T,0]$ with $[T,T,0]M \neq 0$. Pick and fix $m \in
M$ and $T$ such that $m_1 := [T,T,0]m \neq 0$. Then $M = I_0m_1$. 

The elements $[T,T,k]m_1$, $k=0, 1, 2, \ldots $ are not all linearly 
independent. So there is some polynomial $f(X)$ with $f(T, T)m_1 =0$
and it follows that 
$f(S,T)m_1=0$ for all $S\in I_n(0)$, and that
$U_Jm_1=0$ where $J$ is the ideal generated by $f(X)$.  Hence
$U_JM = U_JI_0m_1 \subseteq U_Jm_1=0$. The last part
holds by Remark 2.5.10. 
\qed\enddemo

\bigskip

\subhead 7.2 Indecomposables and cell modules arising from $I_0$
\endsubhead

\definition{7.2.1}Let $\omega(1)$ be the ideal of $D_n$
generated by $(u^n-1)$. Then $u^n$ acts as identity on a module $M$ if
and only if $M$ is a module for the algebra
$D_n/\omega(1)$. We have $\omega(1) \cap I_0 = 0$, so we consider
$I_0$ as a subset (and as an ideal) of $D_n/\omega(1)$. This algebra
is therefore infinite dimensional.
\enddefinition

\definition{7.2.2}Fix an ideal $J$ of $K[X]$, and define an algebra $$
\dnj := \extln/(\omega(1) + U_J)
.$$  This is a finite dimensional algebra, and is the extension of 
$\extln/(\omega(1) + I_0)$ by the ideal $I_0/U_J$. By 7.2.1, every
simple finite dimensional module for $D_n$ with $I_0M \neq 0$ 
is a module for such an algebra $\dnj$. 
Therefore, if we classify the simple modules for these algebras we 
complete the classification of the finite dimensional simple modules for $D_n$.
We will show that the algebra is again cellular, and we will also 
study extensions of its standard modules.
\enddefinition

\proclaim{Lemma 7.2.3}
The algebra $A := \dnj$ is cellular.
\endproclaim

\demo{Proof}
We know that $J_n(1) = D_n/(\omega(1) + I_0)$ is cellular, so 
it suffices to show that 
$I_0/U_J$ has a cell chain. Let $f(X)$ be the monic generator of $J$. Since 
the field is algebraically closed we can factor it into linear factors,
say $$ f(X) = \prod_{i=1}^k (X - r_i)^{t_i}
,$$ where the $r_i$ are taken pairwise distinct and $t_i \geq 1$ for
all $i$.
We can now see that there is a cell basis as for the $q$-Jones algebras 
(see 2.5.1). 
\qed\enddemo

In fact, we can also construct a direct sum decomposition
of $I_0/U_J$ into a direct
sum of indecomposable modules each of which is filtered by standard modules.
This is the same as in \S6.

\definition{7.2.4} 
As in 6.1.1 we can write down
explicit orthogonal primitive idempotents $q_i(X)$ of the ring
$K[X]/\langle f(X) \rangle$, $1 \leq i \leq m$.  Let
$g_i^{(c)}(X) = (X-r_i)^cq_i(X)$, for $0 \leq c \leq a_i$. Define
$Q_i(T)$ to be the span of all $g_c^{(i)}(S, T)$ for $S \in \annn \cap
I_0$ and $0 \leq c \leq a_i$.  Letting $$
U(T) = \bigoplus_{i=1}^m Q_i(T),
$$ we have a direct sum decomposition  $$I_0/U_J = \bigoplus_{T \in
I_n(0)} U(T)
.$$
\enddefinition

\proclaim{Lemma 7.2.5} 
\item{(a)}{The $K$-space $Q_i(T)$ is an $A$-module and it is indecomposable.}
\item{(b)}{There is a filtration of $Q_i(T)$ given by $$
0 \subset V_{N-1} \subset \ldots \subset V_1 \subset V_0 = Q_i(T)
$$ where all the successive quotients are isomorphic.
Moreover, for each $k$ we have $Q_i(T)/V_{N-k} \cong V_k$.}
\endproclaim

\demo{Proof}
This is similar to the proof of Proposition 6.2.1.  We omit details.
\qed\enddemo

\definition{7.2.6}
It follows that up to isomorphism, we obtain one standard module of
$\dnj$ associated to $I_0/U_J$ for each root of $f(X)$, and hence also one simple module for each root.
We label the standard module corresponding
to the root $r_i$ by $W(0, r_i)$. 
Note that $Q_i(T)$ is the projective cover of $W(0, r_i)$. 
\enddefinition

\proclaim{Lemma 7.2.7}
Let $\lambda$ be a weight $(0, r_i)$, and let $\mu = (t, s_j)$ for $0
\leq t$.   Then we have:

\item{(a)}{$\Ext^1(W(\lambda), W(\mu)) = 0$ for $\lambda \neq \mu$;}

\item{(b)}{if $r_i$ has multiplicity $>1$ as a root of $f$ then 
$\Ext^1(W(\lambda), W(\lambda)) = K$,
otherwise $\Ext^1(W(\lambda), W(\lambda)) = 0$.
}
\endproclaim

\demo{Proof}
The proof is the same as that of Theorem 6.3.1 and we omit details.
\qed\enddemo

\example{Example 7.2.8}
Suppose $n=4$, and consider $f = X^N$ for some $N \geq 1$. Then there
is only one standard module up to isomorphism for the ideal
$U_J$.  It is possible for this module not to be simple, as follows.
It is easy to calculate
the Gram matrix of the associated bilinear form as in \cite{{\bf
5}}. One finds that
it has determinant $(\delta^2-2)^2\cdot \delta^4$.  Hence we see that
if $\delta^2 = 2$, the module $W(0,0)$ has a 2-dimensional radical.
\endexample

\vfill\eject
\subhead 7.3 Further properties of modules arising from $I_0$
\endsubhead

\definition{7.3.1} If one considers $\tln$ instead of $D_n$ then the
same techniques are applicable.  One defines algebras
$J_1(n)[J]$ analogous to the algebras $\dnj$.
The only difference is that one has to substitute $[S,T,r(S,T) + 2k]$
for $X^k$ into polynomials. 
\enddefinition

\definition{7.3.2}  
The simple modules constructed for $I_0$ depend only on the root $r$ but not
on the algebra. That is, if $f_1$ and $f_2$ are polynomials which have
both $X-r$ as a factor the the corresponding simple (and standard)
modules of the two algebras with label $(0, r)$ are isomorphic as
modules for $\extln$ (respectively $\tln$). 
This can be seen by using the analogues of 6.1.1 and Lemma 6.2.2.
\enddefinition

\definition{7.3.3}
Let $W = W(0,r)$ be a standard module as above. Then 
$W$ is a module for $\dnj$ (respectively, $J_n(1)[J]$) for any ideal of
$K[X]$ of the form $J = (X-r)^m$ for $m \in \{1, 2, 3, \ldots\}$. 
In this case, Lemma 7.2.5 provides indecomposable modules which have 
filtrations with all quotients isomorphic to $W$, of arbitrary length.
\enddefinition

\definition{7.3.4}
Let $n$ be even. Consider 
$$I_0 + \sp \{ [S,T,w]: w \geq nt \}
.$$
This is a subalgebra of $D_n$ which contains the ideal $I_0$.
We claim that this algebra is also cellular. 
Fix some $q \in K^*$. Then for the quotient modulo $I_0$, which is
$D_n^+$, we take the cell basis as in \S 4. So it suffices to give a cell datum
for $I_0$. Take any linear polynomial $X-r$, and let $J$ be the ideal 
generated by this in $K[X]$. Then we have a chain of ideals
$$\ldots \subset U_{J^k} \subset U_{J^{k-1}} \subset \cdots \subset U_J \subset I_0
$$
which gives a cell basis, as before.


\end